\documentclass [12pt] {article}

\usepackage{titling}
\usepackage{blindtext}
\usepackage {amssymb,latexsym}
\usepackage {amsmath}
\usepackage{graphicx}
\usepackage{amsthm}
\usepackage{relsize}
\usepackage[margin=1in]{geometry}
\usepackage[scr=rsfs,cal=boondox]{mathalfa}
\usepackage{authblk} 
\usepackage{a4wide}
\usepackage[latin1]{inputenc}
\usepackage{fancyhdr}
\usepackage{pgfplots}
\pgfplotsset{compat=1.14}
\usepackage{amscd}
\usepackage{hyperref}
\usepackage{graphicx}
\usepackage{relsize}
\usepackage{newlfont}
\usepackage{amssymb}
\usepackage{youngtab}
\usepackage{young}
\usepackage{mathtools}
\usepackage{tikz}
\usetikzlibrary{shapes, backgrounds} 
\usepackage{tikz-cd}
\usepackage{amsthm}
\usepackage{enumerate,enumitem}

\newtheorem {lemma} {Lemma}
\newtheorem {proposition} [lemma] {Proposition}
\newtheorem {theorem} [lemma] {Theorem}

\newtheorem {definition}[lemma] {Definition}
\theoremstyle{definition}
\newtheorem{remark}[lemma]{Remark}
\newtheorem{example}[lemma]{Example}

\newcommand{\ten}{\otimes}
\newcommand{\id}{\mathrm{id}}
\renewcommand{\b}{\mathfrak{b}}
\newcommand{\cB}{\mathcal{B}}
\newcommand{\cC}{\mathcal{C}}
\newcommand{\cD}{\mathcal{D}}
\newcommand{\cY}{\mathcal{Y}}

\newcommand{\cJ}{\mathcal{J}}

\newcommand{\Com}{\mathrm{Com}}

\begin{document}
\setlength{\parindent}{0pt} 
\title {\textbf{Hopf categories associated to comonoidal functors}}

\author{Andrea Rivezzi} 
\affil{{\sl
{Mathematical Institute of Charles University}}

{\sl Sokolovsk\'a 49/83, 186 75 Prague 8, Czech Republic}}

\begin{titlingpage}

\date{}
\maketitle
\vspace{-1cm}
\begin{abstract}
\noindent
We provide an explicit construction of Hopf categories associated to comonoidal functors,  generalizing \v{S}evera's construction of Hopf monoids through $M$-adapted functors. We discuss the example of the Hopf category whose underlying class is the set of twists of a Lie bialgebra. Finally, we apply the result to the setting of deformed categories.

\end{abstract}

\tableofcontents
\end{titlingpage}

\newpage

%
%
\section*{Introduction}
 \addcontentsline{toc}{section}{Introduction}
 In this paper we discuss Hopf categories, which are situated within the context of of enriched category theory. The latter was introduced by J. B\'enabou \cite{Ben} and J.-M. Maranda \cite{Mar} (see also G.M. Kelly's works \cite{Kelly65,Kelly69,Kelly}) as a generalization of ordinary category theory. 
If $\cC$ is a symmetric monoidal category, one can consider monoid, comonoid, and Hopf monoid objects in $\cC$, by transporting the usual axioms of algebras, coalgebras, and Hopf algebras from $\mathrm{Vect}$ to $\cC$.  The class of the above mentioned objects, together with their class of morphisms form a symmetric monoidal subcategory of $\cC$; hence, one can consider enrichments over such categories. In this paper, we are interested in $\Com(\cC)$-categories,  where $\Com(\cC)$ denotes the category of comonoids of $\cC$. Those are called semi-Hopf $\cC$-categories (or Hopf $\cC$-categories if they admit an antipode), and were introduced by E. Batista, S. Caenepeel, and J. Vercruysse in \cite{BaCaVe} and,  recently, further investigated by several authors  \cite{BFVV,CaeFie,Can,Fie,GroVer, GroVer2}.
 Examples of Hopf categories include the well-known notions of Hopf monoids (where one considers a $\cC$-category with one object) and groupoids (where $\cC = \mathrm{Set}$). \\ \\
In the main result of this paper, see Theorem \ref{main-theorem},  we provide a machinery to construct Hopf categories through the notion of $M$-adapted functor introduced by P. \v{S}evera in \cite{Sev16}.  Namely, if $(F,F^2,F^0)$ is a comonoidal functor and $(M,\Delta,\varepsilon)$ is a comonoid, one says that $F$ is $M$-adapted (and that $M$ is $F$-adapted) if certain morphisms --depending on the comonoidal structure of $F$ and on the comonoid structure of $M$-- are invertible, see Definition \ref{definition-m-adapted-functor}.  The main feature of $M$-adapted functors (which has been used by P. \v{S}evera in order to define a universal quantization functor of Lie bialgebras) is Theorem \ref{theorem-hopf-algebra}, which provides a way to construct, given a cocommutative comonoid $(M,\Delta,\varepsilon)$ and a $M$-adapted functor $F$, a Hopf monoid structure on $F(M \ten M)$. Our idea is to slightly modify this construction in such a way the output is a Hopf category rather than just a Hopf monoid. In order to do so, we consider, given a comonoidal functor $F$, the class $A^0$ of all cocommutative $F$-adapted comonoids.  The comonoidal structure of $F$ allows to define a comonoid object $A_{x,y}$ for each pair of objects $x,y$ in $A^0$.  In order to define the semi-Hopf structure, we then rely on the notion of multiplication along a comonoid introduced in \cite{rivezzithesis}.  The antipode is then defined by means of the braiding of $\cC$. As an application, we show that the set of twists of any Lie bialgebra is the underlying class of a Hopf category. Finally, we extend the main result to the setting of deformed categories.

\subsection*{Acknowledgements}
This paper is inspired by the talk of P. Gro{\ss}kopf at the conference Hopf25 held in Brussels in April 2025.  I wish to thank  M. Bordemann, J. Schnitzer, and T. Weber for numerous enlightening conversations. The author is supported by GA\v{C}R/NCN grant Quantum Geometric Representation Theory and Noncommutative Fibrations 24-11728K and is a member of the \emph{Gruppo Nazionale per le Strutture Algebriche, Geometriche e le loro Applicazioni} (GNSAGA) of the \emph{Istituto Nazionale di Alta Matematica} (INdAM).
This publication is based upon work from \emph{COST Action CaLISTA CA21109} supported by
COST (European Cooperation in Science and Technology). www.cost.eu.

\section{Preliminaries}
\subsection{Symmetric monoidal categories}
We start our discussion by introducing symmetric monoidal categories \cite{MacLane63}, \cite{JS93}.
\begin{definition}
  A braided monoidal category is a datum $ (\mathcal{C},
  \otimes, I, a, \ell, r,\sigma)$, where $\mathcal{C}$ is a category,
  $\otimes : \mathcal{C} \times \mathcal{C} \to \mathcal{C}$ is a
  functor (the tensor product), $I$ is an object in
  $\mathcal{C}$ (the monoidal unit), and
  \begin{equation*}
    \begin{split}
      a&: \ten \circ (\ten \times \id) \Rightarrow \ten \circ (\id \times \ten) \quad \text{  (the associativity constraint)}\\
      \ell&: \ten \circ (I \times \id) \Rightarrow \id \quad \quad \qquad \qquad 
      \text{(the left unit constraint)}\\
      r &: \ten \circ (\id \times I) \Rightarrow \id \qquad \qquad \quad \quad 
      \text{(the right unit constraint)}\\
      \sigma &: \ten \Rightarrow \ten^{\mathrm{op}} \qquad \qquad \quad \quad \quad \qquad
      \text{   (the commutativity constraint, or the braiding)}
    \end{split}
  \end{equation*}
are natural isomorphisms satisfying the following conditions for all objects $X,Y,Z,W$: 
\begin{align*}
(\id_X \ten a_{Y,Z,W}) \circ a_{X,Y \ten Z,W} \circ (a_{X,Y,Z} \ten \id_W) &= a_{X,Y,Z \ten W} \circ a_{X \ten Y,Z,W} \quad \qquad \text{(pentagon axiom)}\\
(\id_X \ten \ell_Y) \circ a_{X,I,Y} &= r_{X} \ten \id_Y \qquad \qquad \qquad \qquad \text{(triangle axiom)}\\
(\id_Y \ten \sigma_{X,Z} )\circ a_{Y,X,Z} \circ (\sigma_{X,Y} \ten \id_Z) &= a_{Y,Z,X} \circ \sigma_{X, Y \ten Z} \circ a_{X,Y,Z} \quad \text{(1st hexagon axiom)}\\
(\sigma_{X,Z} \ten \id_Y)\circ a_{X,Z,Y}^{-1} \circ (\id_X \ten \sigma_{Y,Z}) &= a^{-1}_{Z,X,Y} \circ \sigma_{X \ten Y, Z} \circ a^{-1}_{X,Y,Z}\quad \text{(2nd hexagon axiom)}
\end{align*}  

  If moreover $\sigma^{-1}_{X,Y} =\sigma_{Y,X} $ for all objects $X,Y$
  of $\mathcal{C}$, we say that $\mathcal{C}$ is a symmetric monoidal category.  $\cC$ is said to be strict if the natural isomorphisms $a,\ell,r$ are identities in $\cC$.
\end{definition}
From now on, in view of MacLane's coherence theorem (see e.g.  \cite{MacLane63}, \cite{Kelly64}, \cite{MacLane71}, \cite[XI.5]{Kassel95}) we shall only consider strict symmetric monoidal categories, even though all the results hold in the non-strict setting, and many of them hold in the non-symmetric case. In order to lighten the notation, if $\cC$ and $\cB$ are two symmetric monoidal categories, we shall not use different symbols for their tensor products, monoidal units, and constraints. \\ \\
The usual notion of functor between symmetric monoidal categories is the one of braided monoidal functor. However, we shall need its dual counterpart, i.e. the one of braided \emph{co}monoidal functor, since we shall focus our attention to comonoids.
\begin{definition}
  \label{definition-monoidal-functor}
  Let $\cC$ and $\cB$ be two symmetric monoidal categories. A braided comonoidal functor from $\cC$ to $\cB$ is a triple $(F,F^2,F^0)$, where $F:\cC \rightarrow \cB$ is a functor,   $F^0: F(I_\cC) \to I_\cB$ is a morphism, and $F^2$ is a natural transformation $F^2 : F \circ \ten \Rightarrow \ten \circ (F \times F)$ making commutative the following diagrams for all objects $X,Y,Z$ in $\cC$
\begin{equation*}
\begin{tikzcd}
F(X \ten Y \ten Z) \arrow[d, "{F^2(X \ten Y,Z)}"'] \arrow[rr, "{F^2(X, Y \ten Z)}"] &  & F(X) \ten F(Y \ten Z) \arrow[d, "{\id \ten F^2(Y,Z)}"] &  & F(X \ten Y) \arrow[r, "{F^2(X,Y)}"] \arrow[d, "{F(\sigma_{X,Y})}"'] & F(X) \ten F(Y) \arrow[d, "{\sigma_{F(X),F(Y)}}"] \\
F(X \ten Y) \ten Z \arrow[rr, "{F^2(X,Y) \ten \id}"]                                &  & F(X) \ten F(Y) \ten F(Z)                               &  & F(Y \ten X) \arrow[r, "{F^2(Y,X)}"]                                 & F(Y) \ten F(X)                                  
\end{tikzcd}
\end{equation*}
  and satisfying $(F^0 \ten \id) \circ F^2(I,X) = \id$ and $(\id \ten F^0) \circ F^2(X,I)=\id$.
\end{definition}

\subsection{Hopf monoids in symmetric monoidal categories}
The aim of this section is to define monoids, comonoids, and Hopf monoids in the context of a generic symmetric monoidal category.  It is clear that setting $\cC = \mathrm{Vect}$ one recovers the usual algebraic definitions of algebras, coalgebras and Hopf algebras.
\begin{definition}
Let $\cC$ be a symmetric monoidal category.
\begin{itemize}
\item[(i):] A monoid is a triple $(A, \mu, \eta)$, where $A$ is an object of $\cC$, and $\mu : A \ten A \to A$,  $\eta : I \to A$
are morphisms (called the multiplication and the unit) such that
\begin{align*}
\mu \circ (\mu \ten \id) &= \mu \circ (\id \ten \mu)  \\
\mu \circ (\eta \ten \id) &= \id = \mu \circ (\id \ten \eta) .
\end{align*}
If moreover $\mu \circ \sigma = \mu$ holds, we say that $(A, \mu,\eta)$ is a commutative monoid. 
\item[(ii):] A comonoid is a triple $(M,\Delta, \varepsilon)$, where $M$ is an object of $\cC$, and $\Delta : M \rightarrow M \ten M$, $\varepsilon : M \rightarrow I$ are morphisms (called the comultiplication and the counit) such that
\begin{align*}
(\id \ten \Delta) \circ \Delta&=  (\Delta \ten \id) \circ \Delta  \\
(\varepsilon \ten \id) \circ \Delta &= \id =(\id \ten \varepsilon) \circ \Delta.
\end{align*}
If moreover $\sigma \circ \Delta = \Delta$ holds, we say that $(M,\Delta, \varepsilon)$ is a cocommutative comonoid. 
\item[(iii):] Let $(M, \Delta,\varepsilon)$ and $(M',\Delta',\varepsilon')$ be two comonoids. A morphism $f: M \to M'$ is said to be a morphism of comonoids if $(f \ten f) \circ \Delta = \Delta' \circ f$ and $\varepsilon' \circ f = \varepsilon$.
\end{itemize}
\end{definition}
The following proposition collects some well-known properties of comonoids and comonoidal functors we shall need. 
\begin{proposition}
\label{proposition-properties-comonoids}
Let $\cC$ be a symmetric monoidal category.
\begin{itemize}
\item[(i):] The monoidal unit $I$ is a comonoid with $\Delta = \ell^{-1}_I = r^{-1}_I$ and $\varepsilon = \id$. 
\item[(ii):] If $(M,\Delta,\varepsilon)$ and $(M',\Delta',\varepsilon')$ are two comonoids, then $M \ten M'$ is a comonoid with comultiplication $(\id \ten \sigma \ten \id) \circ (\Delta \ten \Delta')$ and counit $\varepsilon \ten \varepsilon'$.
\item[(iii):] If $f,g$ are two morphisms of comonoids, then $f \ten g$ is a morphism of comonoids.
\item[(iv):] If $(M,\Delta,\varepsilon)$ is a comonoid, then $\varepsilon$ is a morphism of comonoids.
\item[(v):] If $(M,\Delta,\varepsilon)$ is cocommutative comonoid, then $\Delta$ is a morphism of comonoids.
\item[(vi):] If $(M,\Delta,\varepsilon)$ is a comonoid and $(F, F^2,F^0)$ is a comonoidal functor, then $F(M)$ is a comonoid with comultiplication $F^2(M,M)\circ F(\Delta)$ and counit $F^0 \circ F(\varepsilon)$.
\item[(vii):] If $(M,\Delta,\varepsilon)$ and $(M',\Delta',\varepsilon')$ are two comonoids and $(F, F^2,F^0)$ is a comonoidal functor, then $F^2(M,M')$ and $F^0$ are morphisms of comonoids.
\item[(viii):] If $f$ is a morphism of comonoids and $(F, F^2,F^0)$ is a comonoidal functor, then $F(f)$ is a morphism of comonoids.
\end{itemize}
\end{proposition}
Note that statements $(i)-(iii)$ imply that the class of all comonoids forms a symmetric monoidal category. We shall denote it by $\mathrm{Com}(\cC)$. We can thus consider the following 
\begin{definition}
Let $\cC$ be a symmetric monoidal category.
\begin{itemize}
\item[(i):] A bimonoid is a quintuple $(H,\mu,\eta,\Delta,\varepsilon)$, where $(H, \mu,\eta)$ is a monoid, $(H, \Delta,\varepsilon)$ is a comonoid, and $\mu,\eta$ are morphisms of comonoids.
\item[(ii):] A Hopf monoid is a bimonoid with an antipode, that is a morphism $S: H \to H$ making commutative the following diagram
\begin{equation*}
\begin{tikzcd}
                                                                       & H \ten H \arrow[rr, "\id \ten S"] &                      & H \ten H \arrow[rd, "\mu"]  &   \\
H \arrow[ru, "\Delta"] \arrow[rd, "\Delta"'] \arrow[rr, "\varepsilon"] &                                   & I \arrow[rr, "\eta"] &                             & H \\
                                                                       & H \ten H \arrow[rr, "S \ten \id"] &                      & H \ten H \arrow[ru, "\mu"'] &  
\end{tikzcd}
\end{equation*}
\end{itemize}
\end{definition}
\subsection{Enriched categories}
We now present the notion of enrichement over a symmetric monoidal category introcuced by J. B\'enabou \cite{Ben} and J.-M. Maranda \cite{Mar} and further developed by G.M. Kelly \cite{Kelly65,Kelly69,Kelly}. We shall use the same convention as in \cite{BaCaVe} and \cite{GroVer} for the direction of the arrows. 
\begin{definition}
Let $\cC $ be a symmetric monoidal category. A category enriched over $\cC$, or a $\cC$-category, consists of a datum $(A^0,A^1, m,u)$, where:
\begin{itemize}
\item[(i):] $A^0$ is a class, called the underlying class of $(A^0,A^1, m,u)$;
\item[(ii):] $A^1= \{ A_{x,y}\}_{x,y \in A^0}$ is a class of objects in $\cC$ indexed by all pairs in $A^0$;
\item[(iii):] $m = \{m_{x,y,z}  \}_{x,y,z \in A^0}$ is a collection of morphisms $m_{x,y,z} : A_{x,y} \ten A_{y,z} \to A_{x,z}$ in $\cC$ indexed by all triples in $A^0$;
\item[(iv):] $u = \{u_x\}_{x \in A^0}$ is a collection of morphisms $u_x : I \to A_{x,x}$ in $\cC$ indexed by $A^0$
\end{itemize}
such that the diagrams 
\begin{equation}
\label{eq:diagram-semi-hopf-one}
\begin{tikzcd}
{A_{x,y} \ten A_{y,z} \ten A_{z,w}} \arrow[rr, "{m_{x,y,z} \ten \id}"] \arrow[d, "{\id \ten m_{y,z,w}}"'] &  & {A_{x,z} \ten A_{z,w}} \arrow[d, "{m_{x,z,w}}"] \\
{A_{x,y} \ten A_{y,w}} \arrow[rr, "{m_{x,y,w}}"]                                                          &  & {A_{x,w}}                                      
\end{tikzcd}
\end{equation}
and
\begin{equation}
\label{eq:diagram-semi-hopf-two}
\begin{tikzcd}
{A_{x,x} \ten A_{x,y}} \arrow[r, "m_{x,x,y}"]                     & {A_{x,y}} & {A_{x,y} \ten A_{y,y}} \arrow[l, "m_{x,y,y}"']                  \\
{I \ten A_{x,y}} \arrow[ru, "\cong"'] \arrow[u, "u_x \ten \id"] &            & {A_{x,y} \ten I} \arrow[lu, "\cong"] \arrow[u, "\id \ten u_y "']
\end{tikzcd}
\end{equation}
commute for all $x,y,z,w$ in $A^0$.
\end{definition}
\begin{example}
The following are well-known notions of enrichments:
\begin{itemize}
\item[(i):] A $\cC$-category with one object is a monoid in $\cC$.
\item[(ii):] If $\cC= \mathrm{Set}$ then a $\cC$-category is an ordinary category.
\item[(iii):] If $\cC$ is the category of abelian groups, then a $\cC$-category is a preadditive category.
\item[(iv):] If $\cC = \mathrm{Vect}_\mathbb{K}$,  then a $\cC$-category is a $\mathbb{K}$-linear category.
\item[(v):] Another example is given by F.W. Lawvere in \cite{Law}, where he defined pseudometric spaces as categories enriched over the extended real numbers.
\end{itemize}
\end{example}

\subsection{Semi-Hopf categories} 
\label{section-Hopf-cat}
Let $\cC$ be a symmetric monoidal category. Then, since $\mathrm{Com}(\cC)$ is a symmetric monoidal category, we can consider the following
\begin{definition}
A semi-Hopf $\cC-$category is a $\Com(\cC)$-category $(A^0,A^1,m,u)$. If moreover there exists a collection $S=\{S_{x,y}\}_{x,y \in A^0}$ of morphisms $S_{x,y} : A_{x,y} \to A_{y,x}$ in $\cC$ making commutative the diagram
\begin{equation}
\label{eq:Hopf-category}
\begin{tikzcd}
{A_{x,y} \ten A_{x,y}} \arrow[r, "{\id \ten S_{x,y}}"]                                               & {A_{x,y} \ten A_{y,x}} \arrow[r, "{m_{x,y,x}}"] & {A_{x,x}}                            \\
{A_{x,y}} \arrow[rr, "{\varepsilon_{x,y}}"] \arrow[d, "{\Delta_{x,y}}"'] \arrow[u, "{\Delta_{x,y}}"] &                                                 & I \arrow[d, "u_y"] \arrow[u, "u_x"'] \\
{A_{x,y} \ten A_{x,y}} \arrow[r, "{S_{x,y} \ten \id}"]                                               & {A_{y,x} \ten A_{x,y}} \arrow[r, "{m_{y,x,y}}"] & {A_{y,y}}                           
\end{tikzcd}
\end{equation}
we say that $(A^0,A^1,m,u,S)$ is a Hopf $\cC$-category.
\end{definition}
In other words, a semi-Hopf category consists of a class $A^0$ such that:
\begin{itemize}
\item for any $x,y$ in $A^0$, we have a comonoid $(A_{x,y}, \Delta_{x,y}, \varepsilon_{x,y})$ i.e. the following diagrams commute
\begin{equation*}
\begin{tikzcd}
{A_{x,y}} \arrow[r, "{\Delta_{x,y}}"] \arrow[d, "{\Delta_{x,y}}"'] & {A_{x,y} \ten A_{x,y}} \arrow[d, "{\Delta_{x,y} \ten \id}"] &  & {A_{x,y}} & {A_{x,y} \ten A_{x,y}} \arrow[l, "{\id \ten \varepsilon_{x,y}}"'] \arrow[r, "{\varepsilon_{x,y} \ten \id}"] & {A_{x,y}} \\
{A_{x,y} \ten A_{x,y}} \arrow[r, "{\id \ten \Delta_{x,y}}"]        & {A_{x,y} \ten A_{x,y} \ten A_{x,y}}                         &  &           & {A_{x,y}} \arrow[lu, "\cong"] \arrow[ru, "\cong"'] \arrow[u, "{\Delta_{x,y}}"]                                  &          
\end{tikzcd}
\end{equation*}
\item for any $x,y,z$ in $A^0$, we have a morphism $m_{x,y,z}: A_{x,y} \ten A_{y,z} \to A_{x,z}$ satisfying \eqref{eq:diagram-semi-hopf-one} and making commutative the following diagrams:
\begin{equation*}
\begin{tikzcd}
{A_{x,y} \ten A_{y,z}} \arrow[rr, "{\Delta_{x,y} \ten \Delta_{y,z}}"] \arrow[dd, "{m_{x,y,z}}"'] &  & {A_{x,y} \ten A_{x,y} \ten A_{y,z}\ten A_{y,z}} \arrow[d, "\id \ten \sigma \ten \id"]    &  & {A_{x,y} \ten A_{y,z}} \arrow[r, "{\varepsilon_{x,y} \ten \varepsilon_{y,z}}"] \arrow[d, "{m_{x,y,z}}"'] & I \ten I \arrow[d, "\cong"] \\
                                                                                                 &  & {A_{x,y} \ten A_{y,z} \ten A_{x,y} \ten A_{y,z}} \arrow[d, "{m_{x,y,z} \ten m_{x,y,z}}"] &  & {A_{x,z}} \arrow[r, "{\varepsilon_{x,z}}"]                                                               & I                         \\
{A_{x,z}} \arrow[rr, "{\Delta_{x,z}}"]                                                           &  & {A_{x,z} \ten A_{x,z}}                                                                   &  &                                                                                                          &                          
\end{tikzcd}
\end{equation*}
\item for any $x$ in $A^0$, we have a morphism $u_x: I \to A_{x,x}$ satisfying \eqref{eq:diagram-semi-hopf-two} and making commutative the following diagrams: 
\begin{equation*}
\begin{tikzcd}
I \arrow[d, "u_x"'] \arrow[r, "\cong"]  & I \ten I \arrow[d, "u_x \ten u_x"] &  & I \arrow[rr, "{u_{x}}"] \arrow[rd, "\cong"'] &   & {A_{x,x}} \arrow[ld, "{\varepsilon_{x,x}}"] \\
{A_{x,x}} \arrow[r, "{\Delta_{x,x}}"] & {A_{x,x} \ten A_{x,x}}             &  &                                              & I &                                            
\end{tikzcd}
\end{equation*}
\end{itemize}
Semi-Hopf categories were introduced in \cite{BaCaVe} and further investigated in \cite{BFVV,CaeFie,Can,Fie,GroVer, GroVer2}.  We refer to the previous references for a very detailed treatment of them. 
\begin{example}
The following are examples of semi-Hopf categories:
\begin{itemize}
\item[(i):] A semi-Hopf (resp. Hopf) $\cC$-category with one object is a bimonoid (resp. Hopf monoid) in $\cC$.
\item[(ii):] If $\cC = \mathrm{Set}$, then it is shown in \cite[Prop. 2.5]{BaCaVe} that a Hopf $\cC$-category is a groupoid. 
\item[(iii):] A construction of the free and cofree Hopf category is given in \cite[Sect. 4]{GroVer}.
\item[(iv):] It is shown in \cite{GroVer2} that the category of Frobenius algebras is a Hopf category.
\item[(v):] Further examples of semi-Hopf categories, such that the ones of Hopf group algebras and Hopf group coalgebras, are provided in \cite[Sect. 5]{BaCaVe}.
\end{itemize}
\end{example}
\section{The Hopf category of cocommutative $F$-adapted comonoids}
In this section we prove our main result, that is a machinery to associate a Hopf category to any comonoidal functor.  The idea is to generalize \v{S}evera's construction of Hopf monoids (see Theorem \ref{theorem-hopf-algebra}) through the notion of multiplication along a comonoid, see section \ref{section-mult-along-com}. 
\subsection{$M$-adapted functors}
\label{section-M-adapted}
\begin{definition}
\label{definition-m-adapted-functor}
Let $\cC,\cB$ be monoidal categories, $(M, \Delta, \varepsilon)$ be a comonoid in $\cC$ and $(F, F^2,F^0): \cC \to \cB$ be a comonoidal functor. We say that $F$ is $M$-adapted if for any objects $X,Y$ in $\cC$ the following morphisms
\begin{equation*}
\begin{split}
\chi_M &:=F^0 \circ F(\varepsilon): F(M) \to I \\
\gamma_{X,Y}^M &:= F^2(X \ten M, M \ten Y) \circ F\big(\id_X \ten \Delta \ten \id_Y\big) : F(X \ten M \ten Y) \to F(X \ten M) \ten F(M \ten Y)
\end{split}
\end{equation*}
are invertible. Similarly, we say that $M$ is $F$-adapted.
\end{definition}
$M$-adapted functors were introduced  in \cite{Sev16} by P. \v{S}evera in order to construct a universal quantization functor for Lie bialgebras. \v{S}evera's argument is based on the following
\begin{example}
\label{example-lie-bialgebras}
Recall that a Lie bialgebra is a triple $(\mathfrak{b}, [\cdot,\cdot],\delta)$, where $\mathfrak{b}$ is a vector space, and $[\cdot,\cdot]: \b \ten \b \to \b$ and $\delta: \b \to \b \ten \b$ are linear maps, called respectively the Lie bracket and the Lie cobracket, satisfying the following identities:
\begin{equation*}
\begin{split}
[x,y] + [y,x] &=0 \\
[x,[y,z]] + [y,[z,x]] + [z,[x,y]] &=0 \\
\delta(x) + \tau\big(\delta(x)\big) &= 0\\
\mathrm{Alt}\circ (\id_\b \ten \delta) \circ \delta&=0 \\
\delta([x,y]) - x \cdot \delta(y) +  y \cdot \delta(x) &=0
\end{split}
\end{equation*}
where $\tau(x \ten y) = y \ten x$ and $\mathrm{Alt}(x \ten y \ten z) = x \ten y \ten z + y \ten z \ten x + z \ten x \ten y$.
For any Lie bialgebra, one can consider the symmetric monoidal category $\cD \cY (\b)$ of its Drinfeld-Yetter modules. The objects of $\cD \cY(\b)$ are triples $(V, \pi,\pi^*)$, where $V$ is a vector space, and $\pi: \b \ten V \to V$ and $\pi^*: V \to \b \ten V$ are linear maps satisfying
\begin{equation*}
\begin{split}
\pi \circ ([\cdot, \cdot] \ten \id_V) &= \pi \circ (\id_\b \ten \pi) - \pi \circ (\id_\b \ten \pi) \circ (\tau \ten \id_V) \\
(\delta \ten \id_V) \circ \pi^* &= (\tau \ten \id_V) \circ (\id_\b \ten \pi^*) \circ \pi^* -  (\id_\b \ten \pi^*) \circ \pi^* \\
\pi^* \circ \pi = (\id_\b \ten \pi)\circ (\tau &\ten \id) \circ (\id_\b \ten \pi^*) + ([\cdot, \cdot] \ten \id_V)\circ (\id_\b \ten \pi^*) - (\id_\b \ten \pi) \circ (\delta \ten \id_V).
\end{split}
\end{equation*}
The morphisms of $\cD \cY(\b)$ are linear maps which are both a morphism of Lie module and of Lie comodules.
The functor of coinvariants is the braided comonoidal functor
\[ G: \cD \cY(\b) \to \mathrm{Vect}, \quad V \mapsto \frac{V}{\b \cdot V}, \quad f \mapsto G(f)\]
where $G(f)$ is the unique morphism making commutative the diagram
\[
\begin{tikzcd}
V \arrow[r, "f"] \arrow[d, "p_V"'] & W \arrow[d, "p_W"] \\
G(V) \arrow[r, "G(f)"]             & G(W)              
\end{tikzcd}
\]
and $p_V,p_W$ are the canonical projections to the quotient.
One can show that there is a unique Drinfeld-Yetter module structure on the universal enveloping algebra $\mathrm{U}(\b)$ satisfying
\[ \pi(u \ten x) = ux \qquad \text{and} \qquad \pi^*(1)=0.\]
The above structure on $\mathrm{U}(\b)$, together with its standard Hopf algebra structure, is a cocommutative comonoid in $\cD \cY(\b)$.  Moreover, the functor of coinvariants is $\mathrm{U}(\b)$-adapted.
\end{example}
The main feature of $M$-adapted functors is the following theorem, which allows to construct Hopf monoids, see \cite[Sect. 3]{Sev16} for a proof.
\begin{theorem}
\label{theorem-hopf-algebra}
Let $(M, \Delta,\varepsilon)$ be a cocommutative comonoid and $(F, F^2, F^0) : \cC \to \cB$ be a $M$-adapted functor. Then $F(M \ten M)$ is a Hopf monoid, where
\begin{itemize}
\item[(i):] The multiplication is
\begin{equation}
\label{eq:multiplication-hopf-monoid-severa}
 F(\id \ten \varepsilon \ten \id) \circ (\gamma_{M,M}^M)^{-1}.
\end{equation} 
\item[(ii):] The unit is 
\begin{equation}
\label{eq:unit-hopf-monoid-severa}
F(\Delta) \circ (\chi_M)^{-1}.
\end{equation}
\item[(iii):] The comultiplication is 
\begin{equation}
\label{eq:comultiplication-hopf-monoid-severa}
F^2(M \ten M, M \ten M) \circ F(\id \ten \sigma \ten \id) \circ F(\Delta \ten \Delta).
\end{equation}
\item[(iv):] The counit is 
\begin{equation}
\label{eq:counithopf-monoid-severa}
F^0 \circ F(\varepsilon \ten \varepsilon).
\end{equation}
\item[(v):] The antipode is
\begin{equation}
F(\sigma_{M,M}).
\end{equation}
\end{itemize}
\end{theorem}
\subsection{The multiplication along a comonoid}
\label{section-mult-along-com}
In this section we present the notion of multiplication along a $F$-adapted comonoid, which has been introduced by the author in \cite[Sect. 7.2]{rivezzithesis}.
\begin{definition}
Let $\cC, \cB$ be braided monoidal categories, $(M, \Delta,\varepsilon)$ be a comonoid in $\cC$, $(F, F^2, F^0)$ be a $M$-adapted functor from $\cC$ to $\cB$. For any objects $X,Y$ of $\cC$, the multiplication of $X$ and $Y$ along $M$ is the morphism
\begin{equation*}
\mu^M_{X,Y}\coloneqq  F(\id_X \ten \varepsilon \ten \id_Y) \circ (\gamma_{X,Y}^M)^{-1}: F(X \ten M) \ten F(M \ten Y)\to F(X \ten Y) .
\end{equation*}
\end{definition}
Note that the multiplication of the Hopf monoid $F(M \ten M)$ constructed in Theorem \ref{theorem-hopf-algebra} coincides with $\mu_{M,M}^M$. \\
The multiplication along a comonoid satisfies a certain \emph{associative property},  as shown in the following
\begin{theorem}
\label{multiplication-along-a-comonoid-is-associative}
Let $\cC,\cB$ be symmetric monoidal categories, $(M,\Delta_M, \varepsilon_M)$, $(N,\Delta_N,\varepsilon_N)$ be comonoids in $\cC$, and $(F,F^2,F^0): \cC \to \cB$ be a comonoidal functor which is both $M$-adapted and $N$-adapted. Then for any objects $X,Y$ in $\cC$ the following diagram commutes
\begin{equation*}
\begin{tikzcd}
F(X \ten M) \ten F(M \ten N) \ten F(N \ten Y) \arrow[rr, "{\id \ten \mu_{M,Y}^N}"] \arrow[dd, "{\mu_{X,N}^M \ten \id}"'] &  & F(X \ten M) \ten F(M \ten Y) \arrow[dd, "{\mu_{X,Y}^M}"] \\
                                                                                                                         &  &                                                          \\
F(X \ten N) \ten F(N \ten Y) \arrow[rr, "{\mu_{X,Y}^N}"]                                                                 &  & F(X \ten Y)                                             
\end{tikzcd}
\end{equation*}
\end{theorem}
\begin{proof}
The proof follows the same lines of \cite[Rmk. 2]{Sev16}.  Using --for space reasons-- the following notations
\[ F(X) := \overline{X}, \qquad | := \id, \qquad  X \ten Y := XY\] 
we have, recalling that $(F,F^2,F^0)$ is comonoidal, that the diagram
\begin{equation*}
\begin{tikzcd}
\overline{XM} \ \overline{MN} \ \overline{NY}                                                              &  & \overline{XM} \ \overline{MNNY} \arrow[ll, "\id \ten F^2"']                                       &  & \overline{XM} \ \overline{MNY} \arrow[rr, "\id \ten F(| \varepsilon |)"] \arrow[ll, "\id \ten F(| \Delta |)"']                          &  & \overline{XM} \ \overline{MY}                                            \\
                                                                                                           &  &                                                                                                   &  &                                                                                                                                         &  &                                                                          \\
\overline{XMMN} \ \overline{NY} \arrow[uu, "F^2 \ten \id"]                                                 &  & \overline{XMMNNY} \arrow[uu, "F^2"] \arrow[ll, "F^2"']                                            &  & \overline{XMMNY} \arrow[uu, "F^2"] \arrow[ll, "F(|||\Delta|)"'] \arrow[rr, "F(|||\varepsilon|)"]                                        &  & \overline{XMMY} \arrow[uu, "F^2"']                                       \\
                                                                                                           &  &                                                                                                   &  &                                                                                                                                         &  &                                                                          \\
\overline{XMN} \ \overline{NY} \arrow[uu, "F(|\Delta|) \ten \id"] \arrow[dd, "F(|\varepsilon|) \ten \id"'] &  & \overline{XMNNY} \arrow[uu, "F(|\Delta|||)"] \arrow[ll, "F^2"'] \arrow[dd, "F(|\varepsilon|||)"'] &  & \overline{XMNY} \arrow[uu, "F(|\Delta||)"] \arrow[ll, "F(||\Delta|)"'] \arrow[dd, "F(|\varepsilon||)"'] \arrow[rr, "F(||\varepsilon|)"] &  & \overline{XMY} \arrow[uu, "F(|\Delta|)"'] \arrow[dd, "F(|\varepsilon|)"] \\
                                                                                                           &  &                                                                                                   &  &                                                                                                                                         &  &                                                                          \\
\overline{XN} \ \overline{NY}                                                                              &  & \overline{XNNY} \arrow[ll, "F^2"']                                                                &  & \overline{XNY} \arrow[ll, "F(|\Delta|)"'] \arrow[rr, "F(|\varepsilon|)"]                                                                &  & \overline{XY}                                                           
\end{tikzcd}
\end{equation*}
commutes, proving the statement. 
\end{proof}
\subsection{The main result}
\label{section-main-result}
We can now prove the main result of this paper:
\begin{theorem}
\label{main-theorem}
Let $\cC,\cB$ be symmetric monoidal categories, $(F, F^2,F^0): \cC \to \cB$ be a braided comonoidal functor, and $A^0$ be the class of all $F$-adapted cocommutative comonoids.  For any objects $(x,\Delta_x,\varepsilon_x)$, $(y, \Delta_y,\varepsilon_y)$,  $(z, \Delta_z,\varepsilon_z)$ in $A^0$,  set $ A_{x,y} := F(x \ten y)$ and 
\begin{align*}
\Delta_{x,y} &= F^2(x \ten y, x \ten y) \circ F(\id_x \ten \sigma_{x,y} \ten \id_y) \circ F(\Delta_{x} \ten \Delta_y) \\
\varepsilon_{x,y} &= F^0 \circ F(\varepsilon_x \ten \varepsilon_y) \\
m_{x,y,z} &= F(\id_x \ten \varepsilon_y \ten \id_z) \circ (\gamma_{x,z}^y)^{-1}\\
u_x &= F(\Delta_x) \circ (\chi_x)^{-1}\\
S_{x,y} &= F(\sigma_{x,y}).
\end{align*}
Then $(A^0,A^1 = \{ A_{x,y}\},m = \{ m_{x,y,z}\},u = \{u_x\}, S= \{S_{x,y}\})$ is a Hopf category.
\end{theorem}
\begin{proof}
We show that all the axioms of a Hopf category, see section \ref{section-Hopf-cat}, are satisfied.
\begin{itemize}
\item The triple $(A_{x,y}, \Delta_{x,y},\varepsilon_{x,y})$ is a comonoid: since $(x,\Delta_x,\varepsilon_x)$ and $(y, \Delta_y,\varepsilon_y)$ are comonoids, then by part $(ii)$ of Proposition \ref{proposition-properties-comonoids} we have that 
\[\big(x \ten y, (\id_x \ten \sigma_{x,y} \ten \id_y) \circ (\Delta_x \ten \Delta_y), \varepsilon_x \ten \varepsilon_y\big)\] 
is a comonoid.  Since $(F,F^2,F^0)$ is a comomoidal functor, then by part $(vi)$ of Proposition \ref{proposition-properties-comonoids} also 
\[ \big(F(x \ten y), F^2(x \ten y ,x \ten y) \circ F(\id_x \ten \sigma_{x,y} \ten \id_y) \circ F(\Delta_x \ten \Delta_y), F^0 \circ F(\varepsilon_x \ten \varepsilon_y)\big)\]
is a comonoid, which coincides with $(A_{x,y},\Delta_{x,y},\varepsilon_{x,y})$.
\item The morphisms $m_{x,y,z}$ and $u_x$ are morphisms of comonoids: recalling statements $(iii)$, $(iv)$, $(v)$, $(vii)$, $(viii)$ of Proposition \ref{proposition-properties-comonoids},  we have that the morphism
\[m_{x,y,z} = F(\id_x \ten \varepsilon_y \ten \id_z) \circ (\gamma_{x,z}^y)^{-1} = F(\id_x \ten \varepsilon_y \ten \id_z) \circ \big( F^2(x \ten y, y \ten z) \circ F(\id_x \ten \Delta_y \ten \id_z) \big)^{-1} \]
is a morphism of comonoids, since it is a composition of morphisms of comonoids. The same reasoning applies to the morphism
\[ u_x = F(\Delta_x) \circ (\chi_x)^{-1} = F(\Delta_x) \circ \big( F^0 \circ F(\varepsilon_{x,x}) \big)^{-1}. \]
\item The morphism $m_{x,y,z}$ satisfies \eqref{eq:diagram-semi-hopf-one}: it follows by Theorem \ref{multiplication-along-a-comonoid-is-associative} since $m_{x,y,z} = \mu_{x,z}^y$.
\item The morphism $u_x$ satisfies \eqref{eq:diagram-semi-hopf-two}: using the same pictorial notation as in Theorem \ref{multiplication-along-a-comonoid-is-associative}, the claim follows from the commutativity of the following diagram 
\begin{equation*}
\begin{tikzcd}
\overline{xx} \ \overline{xy}                                                                      & \overline{xxxy} \arrow[l, "F^2"']                                                        & \overline{xxy} \arrow[l, "F(|\Delta|)"'] \arrow[r, "F(|\varepsilon|)"]            & \overline{xy} & \overline{xyy} \arrow[l, "F(|\varepsilon|)"'] \arrow[r, "F(|\Delta|)"]           & \overline{xyyy} \arrow[r, "F^2"]                                                        & \overline{xy} \ \overline{yy}                                                                      \\
\overline{x} \ \overline{xy} \arrow[u, "F(\Delta) \ten \id"] \arrow[d, "F(\varepsilon) \ten \id"'] & \overline{xxy} \arrow[l, "F^2"'] \arrow[d, "F(\varepsilon||)"'] \arrow[u, "F(\Delta||)"] & \overline{xy} \arrow[u, "F(\Delta|)"] \arrow[l, "F(\Delta|)"'] \arrow[ru, "\id"'] &               & \overline{xy} \arrow[lu, "\id"] \arrow[u, "F(|\Delta)"'] \arrow[r, "F(|\Delta)"] & \overline{xyy} \arrow[r, "F^2"] \arrow[u, "F(||\Delta)"'] \arrow[d, "F(||\varepsilon)"] & \overline{xy} \ \overline{y} \arrow[u, "\id \ten F(\Delta)"'] \arrow[d, "\id \ten F(\varepsilon)"] \\
\overline{xy} \arrow[d, "F^0 \ten \id"']                                                           & \overline{xy} \arrow[l, "F^2"'] \arrow[ru, "\id"']                                       &                                                                                   &               &                                                                                  & \overline{xy} \arrow[lu, "\id"] \arrow[r, "F^2"]                                        & \overline{xy} \arrow[d, "\id \ten F^0"]                                                            \\
\overline{xy} \arrow[ru, "\id"']                                                                   &                                                                                          &                                                                                   &               &                                                                                  &                                                                                         & \overline{xy} \arrow[lu, "\id"]                                                                   
\end{tikzcd}
\end{equation*}
whose subdiagrams commute since $(F,F^2,F^0)$ is a comonoidal functor and $(x,\Delta_x,\varepsilon_x)$ and $(y, \Delta_y,\varepsilon_y)$ are comonoids.
\item The morphism $S_{x,y}$ satisfies \eqref{eq:Hopf-category}: note first that 
\begin{equation}
\label{eq:proof-antipode-one}
\begin{split}
u_x \circ \varepsilon_{x,y} &= F(\Delta_x) \circ \big( F^0 \circ F(\varepsilon_x)\big)^{-1} \circ F^0 \circ F(\varepsilon_x \ten \varepsilon_y) \\
&= F(\Delta_x) \circ \big( F^0 \circ F(\varepsilon_x)\big)^{-1} \circ F^0 \circ F(\varepsilon_x ) \circ F(\id_x \ten \varepsilon_y) \\
&= F(\Delta_x) \circ F(\id_x \ten \varepsilon_y) \\
&= F(\Delta_x \ten \varepsilon_y)
\end{split}
\end{equation}
and similarly $u_y \circ \varepsilon_{x,y} = F(\varepsilon_x \ten \Delta_y)$. Next, using the naturality of the braiding, the naturality of $F^2$ and the first hexagon axiom, we have that the diagram
\begin{equation*}
\begin{tikzcd}
\overline{xy} \arrow[rr, "F(\Delta \varepsilon)"] \arrow[d, "F(\Delta |)"'] &  & \overline{xx}                                                                \\
\overline{xxy} \arrow[rr, "F(|\sigma)"] \arrow[d, "F(||\Delta)"']               &  & \overline{xyx} \arrow[u, "F(|\varepsilon |)"'] \arrow[dd, "F(|\Delta|)"] \\
\overline{xxyy} \arrow[d, "F(|\sigma|)"']                                         &  &                                                                              \\
\overline{xyxy} \arrow[rr, "F(||\sigma)"] \arrow[d, "F^2"']                       &  & \overline{xyyx} \arrow[d, "F^2"]                                             \\
\overline{xy} \ \overline{xy} \arrow[rr, "\id \ten F(\sigma)"]                    &  & \overline{xy} \ \overline{yx}                                               
\end{tikzcd}
\end{equation*}
commutes. Therefore, we have 
\begin{equation}
\label{eq:proof-antipode-two}
\begin{split}
F(\Delta_x \ten \varepsilon_y) &= \bigg(F(\id_x \ten \varepsilon_y \ten \id_x) \circ \big( F^2(x \ten y,y \ten x) \circ F(\id_x \ten \Delta_y \ten \id_x)\big)^{-1}  \\
& \quad\circ \big(\id_{F(x \ten y)} \ten F(\sigma_{x,y})\big)  \circ F^2(x \ten y, x \ten y) \circ  F(\id_x \ten \sigma_{x,y} \ten \id_y) \circ F(\Delta_x \ten \Delta_y) \bigg)\\
&= m_{x,y,x} \circ (\id \ten S_{x,y}) \circ \Delta_{x,y},
\end{split}
\end{equation}
where we used that the comonoid $(y,\Delta_y,\varepsilon_y)$ is $F$-adapted to invert the morphism $F^2(x \ten y,y \ten x) \circ F(\id \ten \Delta_y \ten \id)$.
Combining equations \eqref{eq:proof-antipode-one} and \eqref{eq:proof-antipode-two} gives \[u_x \circ \varepsilon_{x,y} = m_{x,y,x} \circ (\id \ten S_{x,y}) \circ \Delta_{x,y}\] 
and in a similar way one can show that 
\[u_y \circ \varepsilon_{x,y} = m_{y,x,y} \circ (S_{x,y} \ten \id) \circ \Delta_{x,y}.\]
\end{itemize}
\end{proof}
\begin{remark}
Note that for $x=y=z$ one obtains the Hopf monoid of Theorem \ref{theorem-hopf-algebra}. 
\end{remark}
\subsection{The Hopf category of Lie bialgebra twists}
\label{section-twists}
In this section we apply Theorem \ref{main-theorem} in order to prove that, given a Lie bialgebra, the set of its twists is the underlying class of a Hopf category.  First, recall the following
\begin{definition}
Let $(\b, [\cdot,\cdot],\delta)$ be a Lie bialgebra. A Lie bialgebra twist is an element $j \in \Lambda^2(\b)$ satisfying
\[ \big(\mathrm{Alt} \circ( \delta \ten \id)\big) (j) = [j_{12},j_{23}] + [j_{12},j_{13}] + [j_{13},j_{23}].\]
\end{definition}
We denote the class of all twists of $(\b,[\cdot,\cdot],\delta)$ by $\cJ(\b)$.
The following result is well-known:
\begin{proposition}
\label{proposition-twists}
Let $(\b,[\cdot,\cdot],\delta)$ be a Lie bialgebra and $j \in \cJ(\b)$. Then 
\begin{itemize}
\item[(i):] The triple $(\b_j,[\cdot,\cdot],\delta_j)$ is a Lie bialgebra, where $\b=\b_j$ and
\[ \delta_j(x) = \delta(x) +[1 \ten x+x \ten 1,j]. \]
\item[(ii):] Let $(V,\pi,\pi^*)$ be a Drinfeld-Yetter module of $\b$. Then $(V_j, \pi,\pi^*_j)$ is a Drinfeld-Yetter module for $\b_j$, where $V_j = V$ and 
\[ \pi^*_j(v)= \pi^*(v) + (\id \ten \pi) (j \ten v). \]
\item[(iii):] The functor 
\[ T_j : \cD \cY (\b) \to \cD \cY (\b_j), \qquad V \mapsto V_j , \qquad f \mapsto f\]
is an invertible monoidal functor, and its inverse is $T_{-j}$.
\end{itemize}
\end{proposition}
Following \cite[Rmk. 9]{Sev16}, let $(\b, [\cdot,\cdot],\delta)$ be a Lie bialgebra and consider the Drinfeld-Yetter module structure on the universal enveloping algebra $\mathrm{U}(\b)$ as in Example \ref{example-lie-bialgebras}.  Next, let $j \in \cJ(\b)$ and consider the twisted Lie bialgebra $(\b_j,[\cdot,\cdot],\delta_j)$, together with the Drinfeld-Yetter module $\mathrm{U}(\b_j)$ in $\cD \cY (\b_j)$. By part $(iii)$ of Proposition \ref{proposition-twists} one can consider 
\[ \mathrm{U}_j(\b) := T_{-j}\big( \mathrm{U}(\b_j)\big) \in \cD \cY (\b),\]
which is the unique Drinfeld-Yetter module structure on $\mathrm{U}(\b)$ such that $\pi(x \ten u) = xu$ and $\pi^*(1)=j$. We have the following
\begin{proposition}
Let $(\b,[\cdot,\cdot],\delta)$ be a Lie bialgebra, $j \in \cJ(\b)$,  and consider the functor of coinvariants $G: \cD \cY (\b) \to \mathrm{Vect}$ of Example \ref{example-lie-bialgebras}. Then $\mathrm{U}_j(\b)$ is $G$-adapted. 
\end{proposition}
\begin{proof}
It suffices to observe that the comonoid structure of $\mathrm{U_j(\b)}$ is the same of $\mathrm{U(\b)}$. 
\end{proof}
Therefore, by Theorem \ref{main-theorem} we have that for any Lie bialgebra $(\b, [\cdot,\cdot],\delta)$ there is a Hopf category whose underlying class is $\cJ(\b)$.
\subsection{Deformations of $M$-adapted functors}
In this section we extend Theorem \ref{main-theorem} to the framework of deformed categories. \\
We fix a formal parameter $\hbar$ and a field $\mathbb{K}$ of characteristic zero, and we consider $\mathbb{K}$-linear categories.  Recall the following definitions from \cite{ABSW} and \cite{Sev16}:
\begin{definition}
Let $\cC, \cB$ be braided monoidal categories.
\begin{itemize}
\item[(i)]
We say that $\cC$ is pre-Cartier if there is a natural morphism $t : \ten \Rightarrow \ten$, called the infinitesimal braiding, satisfying
\begin{align*}
t_{X,Y \ten Z} &= t_{X,Y} \ten \id_Z + (\sigma_{X,Y}^{-1} \ten \id_Z) \circ (\id_Y \ten t_{X,Z}) \circ  (\sigma_{X,Y} \ten \id_Z)\\
t_{X \ten Y, Z} &= \id_X \ten t_{Y,Z} + (\id_X \ten \sigma^{-1}_{Y,Z}) \circ (t_{X,Z} \ten \id_Y) \circ  (\id_X \ten \sigma_{Y,Z}) .
\end{align*}
\item[(ii):] A comonoid $(M,\Delta,\varepsilon)$ in a pre-Cartier category is said to be infinitesimally cocommutative if $t \circ \Delta = \sigma \circ \Delta = \Delta$.
\item[(iii):] A braided comonoidal functor $(F,F^2,F^0): \cC \to \cB$ between pre-Cartier categories is said to be infinitesimally braided if the following diagram commutes
\begin{equation*}
\begin{tikzcd}
F(X  \ten Y) \arrow[r, "{F(t_{X,Y})}"] \arrow[d, "{F^2(X,Y)}"'] & F(X  \ten Y) \arrow[d, "{F^2(X,Y)}"] \\
F(X ) \ten F(Y) \arrow[r, "{t_{F(X),F(Y)}}"]                    & F(X ) \ten F(Y)                     
\end{tikzcd}
\end{equation*}
\end{itemize}
\end{definition}
The existence of an infnitesimal braiding in a braided monoidal category allows to consider deformations \emph{\`a la Cartier},  see \cite{ERSW} for precise definitions.  Recall the following result, see \cite[Thm. 4.11]{ERSW} for a proof.
\begin{theorem}
Let $\cC$ be a pre-Cartier category satisfying
\begin{equation}
\label{eq:commutation-relation}
 [t \ten \id , \id \ten t]=0 .
\end{equation}
Then there exists a deformation of $\cC$, which we denote by $\hat{\cC}$, whose braiding is given by 
\begin{equation}
\label{eq:braiding-deformed}
 \hat{\sigma }_{X,Y} = \tilde{\sigma}_{X,Y} \circ e^{\hbar t_{X,Y}}
\end{equation}
where $\tilde{\sigma}$ denotes the $\hbar$-linear extension of $\sigma$.
\end{theorem}
An easy computation shows that the braiding \eqref{eq:braiding-deformed} is symmetric if and only if
\begin{equation}
\label{eq:condition-for-symmetric}
t_{Y,X} \circ \sigma_{X,Y} = - \sigma_{X,Y} \circ t_{X,Y}.
\end{equation}
Next, consider two pre-Cartier categories $\cC,\cB$,  an infinitesimally cocommutative comonoid $(M,\Delta,\varepsilon)$ in $\cC$, and a $M$-adapted infinitesimally braided monoidal functor $(F, F^2,F^0): \cC \to \cB$. Then \v{S}evera showed (see \cite[Prop. 2]{Sev16}) that $F$ extends to a deformed comonoidal functor $(\hat{F}, \hat{F}^2, \hat{F}^0): \hat{\cC} \to \hat{\cB}$, which turns out to be again $M$-adapted.  Therefore,  we can apply Theorem \ref{main-theorem} to the class of all infinitesimally cocommutative comonoids adapted to $F$ and form a new Hopf category.  In other words, we have the following
\begin{theorem}
\label{main-theorem-deformed}
Let $\cC,\cB$ be symmetric pre-Cartier categories satisfying conditions \eqref{eq:commutation-relation} and \eqref{eq:condition-for-symmetric}, $(F, F^2,F^0): \cC \to \cB$ be an infinitesimally braided comonoidal functor, and $A^0$ be the class of all $F$-adapted infinitesimally cocommutative comonoids.  For any objects $(x,\Delta_x,\varepsilon_x)$, $(y, \Delta_y,\varepsilon_y)$,  $(z, \Delta_z,\varepsilon_z)$ in $A^0$,  set $ A_{x,y} := \hat{F}(x \ten y)$ and 
\begin{align*}
\Delta_{x,y} &= \hat{F}^2(x \ten y, x \ten y) \circ \hat{F}(\id_x \ten \hat{\sigma}_{x,y} \ten \id_y) \circ \hat{F}(\Delta_{x} \ten \Delta_y) \\
\varepsilon_{x,y} &= \hat{F}^0 \circ \hat{F}(\varepsilon_x \ten \varepsilon_y) \\
m_{x,y,z} &= \hat{F}(\id_x \ten \varepsilon_y \ten \id_z) \circ (\gamma_{x,z}^y)^{-1}\\
u_x &= \hat{F}(\Delta_x) \circ (\chi_x)^{-1}\\
S_{x,y} &= \hat{F}(\hat{\sigma}_{x,y}).
\end{align*}
Then $(A^0,A^1 = \{ A_{x,y}\},m = \{ m_{x,y,z}\},u = \{u_x\},S = \{ S_{x,y}\})$ is a Hopf category.
\end{theorem}

%
%
%


\begin{thebibliography}{99}

\bibitem{ABSW}
A. Ardizzoni, L. Bottegoni, A. Sciandra, T. Weber:
\emph{Infinitesimal braidings and pre-Cartier bialgebras}.
Communications in Contemporary Mathematics 27, 5 (2025) 2450029.

\bibitem{BaCaVe}
E. Batista, S. Caenepeel, J. Vercruysse:
\emph{Hopf categories}.
Algebras and representation theory 19 (2016): 1173-1216.

\bibitem{Ben}
J. B\'enabou:
\emph{Cat\'egories relatives}.
Comptes Rendus de l'Acad\'emie des Sciences 260 (1965), pp. 3824-3827.

\bibitem{BFVV}
M. Buckley, T. Fieremans, C. Vasilakopoulou, J. Vercruysse:
\emph{A Larson-Sweedler theorem for Hopf $\mathcal{V}$-categories.} Advances in Mathematics, 376, 107456 (2021).

\bibitem{CaeFie}
S. Caenepeel, T. Fieremans:
\emph{Descent and Galois theory for Hopf categories.}
Journal of Algebra and its Applications 17.07 (2018): 1850120.

\bibitem{Can}
C. R. Canlubo:
\emph{Hopf algebroids, Hopf categories and their Galois theories.} Preprint arXiv:1612.06317 (2016).

\bibitem{ERSW}
C. Esposito, A. Rivezzi, J. Schnitzer, T. Weber:
\emph{Quantization of infinitesimal braidings and pre-Cartier quasi-bialgebras}.
Preprint arXiv:2505.17729 (2025).

\bibitem{Fie}
T. Fieremans:
\emph{Hopf and Frobenius $\mathcal{V}$-categories.} (Ph.D. thesis), https://cris.vub.be/ws/portalfiles/portal/65723338/PhDTimmyFieremans.pdf (2019).

\bibitem{GroVer}
P. Gro{\ss}kopf, J. Vercruysse:
\emph{Free and co-free constructions for Hopf categories}.
Journal of Pure and Applied Algebra 228.10 (2024): 107704.

\bibitem{GroVer2}
P. Gro{\ss}kopf, J. Vercruysse:
\emph{The Hopf category of Frobenius algebras}.
Preprint arXiv:2406.18499 (2024).

\bibitem{JS93}
A. Joyal, R. Street:
\emph{Braided tensor categories}. 
Advances in Mathematics 102 (1993) 20-78.

\bibitem{Kassel95}
C. Kassel: \emph{Quantum groups}. Graduate Texts in Mathematics, 155. Springer-Verlag, New York, 1995.

\bibitem{Kelly64}
G.M. Kelly:
\emph{On MacLane's conditions for coherence of natural associativities, commutativities, etc.}
Journal of Algebra, 1(4), 397-402. (1964).

\bibitem{Kelly65}
G.M. Kelly:
\emph{Tensor products in categories}.
 Journal of Algebra 2(1965), 15-37.
 
\bibitem{Kelly69}
G.M. Kelly:
\emph{Adjunction for enriched categories}.
Lecture Notes in Mathematics 420(1969), 166-177.

\bibitem{Kelly}
G. M. Kelly:
\emph{Basic concepts of enriched category theory}.
Vol. 64. CUP Archive, 1982.

\bibitem{Law}
F.W. Lawvere:
\emph{Metric spaces, generalized logic, and closed categories.} Rendiconti del seminario matematico e fisico di Milano 43 (1973): 135-166.

\bibitem{MacLane63} 
S. Mac Lane: 
\emph{Natural Associativity and Commutativity}.  
Rice University Studies 49 (1963) 28-46.

\bibitem{MacLane71}
S. Mac Lane:
\emph{Categories for the Working Mathematician}. volume 5 of Graduate Texts in Math. Springer-Verlag, New York, 1971.

\bibitem{Mar}
J.-M. Maranda:
\emph{Formal categories}.
Canadian Journal of Mathematics 17 (1965) 758-801.

\bibitem{rivezzithesis}
A. Rivezzi:
\emph{Universal constructions arising from quantization of Lie bialgebras}.
 (Ph.D. thesis), www.boa.unimib.it/handle/10281/475779 (2024).

\bibitem{Sev16} P. \v{S}evera: 
\emph{Quantization of Lie bialgebras revisited}.
Selecta Mathematica. 22 (2016), 1563--1581.


\end{thebibliography}
\end{document}